\documentclass[12pt]{article}
\usepackage{amsthm,amstext,amscd,latexsym,amsmath,amsfonts,amssymb}
\unitlength=1mm\linethickness{0.6pt}
\newcommand{\bk}{\Bbbk}
\newcommand{\Z}{\mathbb{Z}}

\newcommand{\N}{\mathbb{N}}
\newcommand{\ie}{\textit{i.e.}\,}

\newcommand{\lin}{\text{lin}}

\newcommand{\id}{\text{id}}

    \theoremstyle{definition}
\newtheorem{def1}{Definition}[section]
\newtheorem{def2}[def1]{Definition}
\newtheorem{def3}[def1]{Definition}
\newtheorem{defsc}[def1]{Definition}
\newtheorem{com1}[def1]{Comment}
\newtheorem{com2}[def1]{Terminology}
\newtheorem{Lemma0}[def1]{Lemma}
\newtheorem{def4}[def1]{Definition}
\newtheorem{def4a}[def1]{Definition}
\newtheorem{opp}[def1]{Definition}
\newtheorem{exm2}[def1]{Example}
\newtheorem{exm4}[def1]{Example}
\newtheorem{exm4a}[def1]{Example}
\newtheorem{exm4b}[def1]{Example}
\newtheorem{def3a}[def1]{Definition}
\newtheorem{def5}{Definition}[section]
\newtheorem{def6}[def5]{Definition}
\newtheorem{def7}[def5]{Definition}
\newtheorem{def8}[def5]{Definition}
\newtheorem{prob1}{Problem}[section]

\newtheorem{exm}[prob1]{Example}
\newtheorem{dfn}[prob1]{Definition}

\begin{document}\title{n-Categories Admissible by
n-graph\thanks{Submitted to `Miscellanea Algebraica', Waldemar
Korczy\'nski and Adam Obtu{\l}owicz, E\-di\-tors, Akademia
\'Swi{\c e}tokrzyska, Kielce, Poland.\hfil\break
Supported by el Consejo Nacional de Ciencia y Tecnolog\'{\i}a (CONACyT)
de M\'exico, proyecto \# 27670 E (1999-2000), and by UNAM, DGAPA,
Programa de Apoyo a Proyectos de Investigacion e Innovacion Tecnologica,
proyecto \# IN 109599.}}
\author{W{\l}adys{\l}aw Marcinek\footnotemark[1]\quad and\quad Zbigniew
Oziewicz\footnotemark[2]}\date{Submitted January 3, 2001}\maketitle
\footnotetext[1]{Institute of Theoretical Physics, University of
Wroclaw, Plac Maksa Borna 9, 50-204 Wroclaw, Poland,
wmar@ift.uni.wroc.pl}\footnotetext[2]{Universidad Nacional Aut\'onoma
de M\'exico, Facultad de Estudios Superiores, Apartado Postal 25, 54700
Cuautitl\'an Izcalli, Estado de M\'exico,\\oziewicz@servidor.unam.mx;
and  Uniwersytet Wroc{\l}awski, Instytut Fizyki Teoretycznej,
Pl 50204 Wroc{\l}aw, Poland.\hfill\break A member of Sistema Nacional de
Investigadores, M\'exico, No de expediente 15337.}
\begin{abstract}We consider examples of the n-categories having the
same n-graph.\end{abstract}\tableofcontents\hyphenation{o-pe-ra-tion}

\section{n-category as n-graph with a composition} We define an
n-category as an n-graph equipped with the partial compositions of
cells, or even with the graftings of cells in case a family of 0-cells
is a strict monoid [Oziewicz \& V\'azquez Couti\~no 1999]. Therefore
here n-graph $G$ is a carrier of n-category $(G,c_G)$ and there is a
forgetful functor from n-category to his n-graph, $(G,c_G)\mapsto G.$

\begin{def1}[Graph] An $\infty$-graph $G$ is an infinite sequence of
families of $i$-cells $\{G_i,i\in\Z\}$ a sequence of sourcees
$\{s_G\equiv s_i\}$ and targets $\{t_G\equiv t_i\}$ surjections, and the
family of sections $\{i_G\equiv i_i\},$

$$\begin{picture}(82,20) \put(0,10){\makebox(0,0)[cc]{$\cdots$}}
\put(5,10){\makebox(0,0)[cc]{$\textit{\large G}_0$}}
\put(29,10){\makebox(0,0)[cc]{$\textit{\large G}_1$}}
\put(53,10){\makebox(0,0)[cc]{$\textit{\large G}_2$}}
\put(77,10){\makebox(0,0)[cc]{$\cdots$}} \put(12,10){\vector(1,0){10}}
\put(11,12){\vector(-2,-1){0.2}}
\bezier{56}(23,12)(17,18)(11,12)\put(11,8){\vector(-2,1){0.2}}
\bezier{56}(23,8)(17,4)(11,8)\put(17,5){\makebox(0,0)[ct]{$t_0$}}
\put(17,16){\makebox(0,0)[cb]{$s_0$}} \put(17,11){\makebox(0,0)[cb]{$i_0$}}
\put(36,10){\vector(1,0){10}} \put(60,10){\vector(1,0){10}}
\put(35,12){\vector(-2,-1){0.2}}
\bezier{56}(47,12)(41,18)(35,12)\put(59,12){\vector(-2,-1){0.2}}
\bezier{56}(71,12)(65,18)(59,12)\put(35,8){\vector(-2,1){0.2}}
\bezier{56}(47,8)(41,4)(35,8)\put(59,8){\vector(-2,1){0.2}}
\bezier{56}(71,8)(65,4)(59,8)\put(41,5){\makebox(0,0)[ct]{$t_1$}}
\put(65,5){\makebox(0,0)[ct]{$t_2$}} \put(41,16){\makebox(0,0)[cb]{$s_1$}}
\put(65,16){\makebox(0,0)[cb]{$s_2$}} \put(41,11){\makebox(0,0)[cb]{$i_1$}}
\put(65,11){\makebox(0,0)[cb]{$i_2$}}\end{picture}$$

Let $\text{type}\,x\equiv(sx,tx).$ By definition all 0-cells are of the same
type,
there are no more then two $(-1)$-cells and at least must be one
$(-1)$-cell,
$G_{-1}\equiv\{s_{-1}G_0,t_{-1}G_0\}.$ A bimap $G_i:G_i\times G_i\rightarrow
2^{G_{i+1}}$ is defined as follows. If $x,y\in G_i\times G_i,$ then
$$G_i(x,y)\equiv\{z\in G_{i+1},\text{type}\,z=(x,y)\}\subset G_{i+1}.$$

The following conditions must hold.\begin{description}
\item[\qquad Type.] If $\text{type}\,x\neq\text{type}\, y,$ then
$G_i(x,y)=\emptyset.$
\item[\qquad Section.] $s_i\circ i_i=\id_{G_i}=t_i\circ i_i.$
\end{description}\end{def1}

\begin{def2}[$n$-graph] If $\forall\,i\in\N,$ $G_{i+n}\simeq G_n,$ then
an $\infty$-graph $G$ is said to be $n$-graph. Hence $s_{i+n}=t_{i+n}$
and $i_{n+i}\circ s_{n+i}=\id_{G_{n+i+1}}.$\end{def2}

In the terminology of Street [1998, p. 100] what we call $n$-graph $G$
is said to be $n$-sceletal globular object, however we do not adopt this
terminology in what follows. For example we do not wish to exclude from
the further consideration more then two surjections
$\{s,t,\ldots\}.$

\begin{defsc}[Sceletal graph] If $\forall\,x,y,$ such that
$\text{type}\,x=\text{type}\,y$ and $\forall\,i\in\N,$
$|G_i(x,y)|=1,$ then a graph $G$ is said to be sceletal.\end{defsc}

\begin{def3}[n-category] An $n$-graph $G$ with a operation
$c_G,$ is said to be an $n$-category,
$$c_G\equiv\left\{\begin{CD}G_i(x,y)\times G_i(y,z)@>{c_ixyz}>>
G_i(x,z)\end{CD}\right\}.$$
In case $c_ixyz$ needs not to be global $\forall\,x,y,z$ then
$n$-category is said to be partial.\end{def3}

Steiner [1998] introduced `genuine' and partial $n$-category as a
$n$-category generated by a `complex'. Therefore the Steiner's complex
is generating a unique $n$-graph admitting a unique $n$-category (either
partial or global. In our definition an $n$-graph do not need admit a
unique $n$-category (either global or partial).

\begin{com1} An n-category $G$ is said to be \textit{strict} if $c_G$
is associative, \ie if a theory is equational. Otherwise an
$n$-category is said to be \textit{weak}. Gray [1976] investigated
tensor pro\-duct of 2-categories. Simpson [1998] studied strict
n-category, 3-grupoids and homotopy types.

A week 2-category is said to be a bicategory [B\'enabou 1967]. A week
3-category is said to be a tricategory [Gordon, Power \& Street 1995;
Baez \& Dolan 1996; Leinster 1998]. Every tricategory is triequivalent
to a category of bicategories equipped with the strong product of Gray
[1976]. Baez and Dolan [1996], Baez [1997] and Leinster [1998],
considered a weak n-category.

Batanin introduced a globular category [Batanin 1998, Street 1998].
\end{com1}

\begin{com2} 0-cell is also called an object, 0-morphism, node; 1-cell
- a morphism, 1-morphism, an arrow, a functor (Definition 3.2). A
2-cell is called a 2-morphism, morphism of morphism, natural
transformation (Definition 3.3), etc.\end{com2}

\begin{Lemma0} A sceletal $n$-graph admits exactly one $n$-category.
\end{Lemma0}

\begin{def4}[Operad = Sketch with graftings] B\'enabou noticed in 1967
that $G_0$ is a monoid iff there is a single $(-1)$-cell only. Then a
n-category $G$ is said to be a monoidal. If $G_0$ is a
\textit{free} monoid on k-generators, then a category $G$ is said to be
a k-sorted sketch. In this case a partial operation $c_G$ from Definition
1.4 generalize to grafting $G_1\times G_1\rightarrow 2^{G_1}$ [Oziewicz
\& V\'azquez 1999].\end{def4}
The above observation by B\'enabou suggest to encode further category
structure in a hidden negative $\Z$-grade `tail' of a graph $G$ as
follows.
\begin{def4a}[Structure tail] A finite sequence of negative $\Z$-grade
cells $\{G_k,-\infty<k\leq-1\}$ with $|G_k|\leq 2,$ is said to be a
structure tail of the category $G.$\end{def4a}

\begin{opp}[Opposite category] Let $C$ be $n$-category. Then $\forall\,
i\geq 1,$ $i$-opposite n-category $C^{opp}_i$ is defined by `reversing'
i-cells, $s^{opp}_i=t_i,$ $t^{opp}_i=s_i,$ and unchanged j-cells for
$j\neq i,$ $s^{opp}_{j\neq i}=s_j,\;t^{opp}_{j\neq i}=t_j.$\end{opp}

\begin{exm2} A collection $G_i(x,y)\subset G_{i+1}$ is a 1-category,
\begin{align*}(G_i(x,y))_0&\equiv G_i(x,y),\\(G_i(x,y))_1&\equiv
\{[G_i(x,y)](\alpha,\beta)\subset G_{i+2}:\quad\forall\;
\alpha,\beta\in G_i(x,y)\}.\end{align*}\end{exm2}

\begin{exm4}Let $\bk$ be a field. Let every 0-cell be a $\bk$-space and
thus $G_0$ be unital $\N$-algebra (unital rig) with $\oplus$ and
$\otimes\equiv\otimes_\bk.$ Let $G_0(V,W)$ be a $\bk$-space
$\lin_\bk(V,W)$ and $G_0(V,W)(f,g)\subset G_1(f,g)\subset G_2$ with
$|G_1(f,g)|=1.$ Then $G$ is a 2-category.\end{exm4}

\begin{exm4a}[Wolff 1974] Enriched Wolff's graph = 2-graph.\end{exm4a}

\begin{exm4b}n-Grupoid is an n-category with an inversion [Simpson 1998].
\end{exm4b}

\begin{def3a}[Cocategory] An n-graph $G$ with a binary cooperation
$\triangle_G$ is said to be an n-cocategory,
$$\triangle_G\equiv\left\{\begin{CD}G_i(x,y)@>{\triangle_ixy}>>
G_i(x,\ldots)\times G_i(\ldots,y)\end{CD}\right\}.$$\end{def3a}

\section{Graph of graphs and category of categories} In this Section
$G$ is $\infty$- or $n$-graph such that $G_0$ is a collection of
k-graphs, so $G$ is a `$n$-graph of $k$-graphs'.

\begin{def5}[Covariant vs contra] Let $E,F\in G_0$ be two k-graphs. An
$k$-graph (absolutely covariant) morphism $f\in G_0(E,F)\subset G_1$ is
a sequence of mappings $\{f_i\}$ intertwining $s_E,t_E,i_E$ with
$s_F,t_F,i_F,$
$$\begin{CD}E_0@<{s_0,t_0}<<E_1@<{s_1,t_1}<<E_2@<{s_2,t_2}<<\ldots,\\
@V{f_0}VV@V{f_1}VV@V{f_2}VV@V{f}VV\\
F_0@<{s_0,t_0}<<F_1@<{s_1,t_1}<<F_2@<{s_2,t_2}<<\ldots.\end{CD}$$
A 1-cell $f\in G_0(E,F^{opp}_i)$ or $\in G_0(E^{opp}_i,F)$ is said to be
$i$-contravariant.\end{def5}
Weakening: if the above intertwining (= a commutativity $f\circ s=s\circ
f$) is up to 2-cell $\alpha:f\circ s\rightarrow s\circ f,$ then $f$ is
said to be pesudo-morphism or $\alpha$-morphism, or weekened morphism.

\begin{def6}[Functor] Let $(E,c_E)$ and $(F,c_F)$ be $k$- or
$\infty$-categories. An $k$-graph morphism $f\in G_0(E,F)$ is said to
be an $k$-category morphism $f\in G_0(c_E,c_F)$ if $f$ intertwine $c_E$
with $c_F,$
$$\begin{CD}E_i(x,y)\times E_i(y,z)@>{c_E}>>E_i(x,z)\\
@V{f_{i+1}\times f_{i+1}}VV@V{f_{i+1}}VV\\
F_i(f_ix,f_iy)\times F_i(f_ix,f_iz)@>{c_F}>>F_i(f_ix,f_iz)
\end{CD}$$\end{def6}

\begin{def7}[Natural transformation, Eilenberg \& Mac Lane 1945] Let
$f,g\in G_0(c_E,c_F).$ A k-tran\-s\-for\-ma\-tion $t\in G_1(f,g),$ is a
collection of $(i+1)$-cells $tx\in F_i(f_ix,g_ix),$ intertwining $f$
with $g,$ $$\begin{CD}E_i(x,y)@>{f_{i+1}}>>F_i(f_ix,f_iy)\\
@V{g_{i+1}}VV@V{(t_y)_*}VV\\
F_i(g_ix,g_iy)@>{(t_x)^{\ast}}>>F_i(f_ix,g_iy)\end{CD}$$
For $\forall\,a\in E_i(x,y),\;(g_{i+1}a)\circ(tx)=(ty)\circ(f_{i+1}a)\;
\in F_{i+1}.$\end{def7}

In the above Definitions 2.1-2.3 we used the (labelled) commutative
dia\-grams as invented by Witold Hurewicz in 1941. In the next Definition
we prefer to use the `graphical calculus', \ie the relations (also known
as `identities'), in an operad of graphs. This calculus is has been
invented among other by Yetter [1990], by Joyal \& Street [1993], and
can be also coined as the Yetter-Joyal-Street \textit{string} calculus.
Calculus of n-graphs is an alternative language for commutative
diagrams. Hurewicz associated to i-cell from $G_i$ the `i-dimensional
geometrical arrow', an oriented solid, \ie 0-dimensional (labelled)
vertex, a point, as a 0-cell $\in G_0,$ a
1-dimensional directed arrow as a 1-cell from $G_1,$ a 2-dimensional
oriented surface as a 2-cell from $G_2,$ etc. The graphical calculus for
an n-graph $G$ consists in representing n-cell as a 0-dimensional
vertex, then (n-1)-cell is represented as a 1-dimensional directed line,
an arrow, etc, and finally 0-cell from $G_0$ is represented as an
oriented n-dimensional volume.

Any segment of two consequent families of n-graph $G,$
$G_i\leftarrow G_{i+1},$ can be represented in the operad of graphs as
follows: if one-dimensional directed edge (an arrow) represents an
i-cell, then zero-dimensional node must represents an (i+1)-cell, \ie an
(i+1)-morphism of i-cells.

Any segment of three families of n-graph $G$ can be represented in these
alternative languages as follows
$$\text{\begin{tabular}{c}The Hurewicz's\\CDs\\\\\\\\The graphical\\
calculus\end{tabular}}\quad\qquad\begin{CD}0@>{+1}>>1@>{+1}>>2@>{+1}>>\\
@A{\text{dim}}AA@A{\text{dim}}AA@A{\text{dim}}AA\\
G_i@<{s_i,t_i}<<G_{i+1}@<{s_{i+1},t_{i+1}}<<G_{i+2}@<<<\\
@V{\text{dim}}VV@V{\text{dim}}VV@V{\text{dim}}VV\\2@<{+1}<<1@<{+1}<<0
\end{CD}$$
All graphs in an operad are implicitely directed from the top (top is an
input) to the bottom (output). Dashed boxes (from infinity) are for
convenience. For two families segments all nodes are
inner, there are no outer nodes. For three families segments all nodes
and arrows are inner, there are neither outer arrows, nor outer nodes
[Oziewicz \& V\'azquez Couti\~no 1999].

Because in our convention all graphs are implicitely directed from the
top to the bottom, therefore what is known in the commutative diagrams
of opetopes [Baez and Dolan 1997], globes [Batanin 1998], and so on, as
the vertical composition, in an operad of graph is the horizontal
composition rather and vice versa. The following Figures illustrate
this correspondence.

\begin{center}\begin{picture}(100,60)
\put(5,10){\dashbox{5}(30,40)[cc]{}} \put(76,37){\line(1,-1){4}}
\put(80,33){\line(1,1){4}} \put(84,17){\line(-1,-1){4}}
\put(80,13){\line(-1,1){4}}
\put(5,15){\line(1,1){30}} \put(15,25){\circle*{1.5}}
\put(25,35){\circle*{1.5}}
\put(50,30){\vector(-1,0){5}} \put(50,30){\vector(1,0){5}}
\put(65,30){\circle*{1.5}}
\put(95,30){\circle*{1.5}} \put(70,30){\vector(1,0){20}}
\put(93,35){\vector(1,-2){0.2}}
\bezier{260}(67,35)(80,65)(93,35)
\put(93,25){\vector(1,2){0.2}}
\bezier{260}(67,25)(80,-5)(93,25)
\put(79,45){\line(0,-1){10}} \put(81,35){\line(0,1){10}}
\put(81,25){\line(0,-1){10}}
\put(79,15){\line(0,1){10}} \put(10,45){\makebox(0,0)[cc]{{\large x}}}
\put(30,15){\makebox(0,0)[cc]{{\large y}}}
\put(62,30){\makebox(0,0)[cc]{{\large x}}}
\put(97,30){\makebox(0,0)[lc]{{\large y}}}
\put(20,51){\makebox(0,0)[cb]{Operad of
graphs}} \put(80,51){\makebox(0,0)[cb]{Globes \& opetopes}}
\put(20,9){\makebox(0,0)[ct]{horizontal}}
\put(80,9){\makebox(0,0)[ct]{vertical}}
\put(29,41){\makebox(0,0)[rb]{{\large f}}}
\put(91,41){\makebox(0,0)[lb]{{\large
f}}}\end{picture}\end{center}

\begin{center}\begin{picture}(120,43)
\put(5,10){\dashbox{5}(20,30)[cc]{}} \put(5,10){\line(1,1){20}}
\put(5,20){\line(1,1){20}} \put(15,20){\circle*{1.5}}
\put(15,30){\circle*{1.5}}
\put(40,25){\vector(-1,0){5}} \put(40,25){\vector(1,0){5}}
\put(55,25){\circle*{1.5}}
\put(85,25){\circle*{1.5}} \put(115,25){\circle*{1.5}}
\put(83,27){\vector(3,-4){0.2}}
\bezier{164}(57,27)(70,43)(83,27)
\put(83,23){\vector(3,4){0.2}} \bezier{164}(57,23)(70,7)(83,23)
\put(69,30){\line(0,-1){10}} \put(71,20){\line(0,1){10}}
\put(66,22){\line(1,-1){4}}
\put(70,18){\line(1,1){4}} \put(113,27){\vector(3,-4){0.2}}
\bezier{164}(87,27)(100,43)(113,27) \put(113,23){\vector(3,4){0.2}}
\bezier{164}(87,23)(100,7)(113,23) \put(99,30){\line(0,-1){10}}
\put(101,20){\line(0,1){10}} \put(96,22){\line(1,-1){4}}
\put(100,18){\line(1,1){4}}
\put(10,35){\makebox(0,0)[cc]{{\large x}}}
\put(10,20){\makebox(0,0)[cc]{{\large y}}}
\put(20,15){\makebox(0,0)[cc]{{\large z}}}
\put(55,26.33){\makebox(0,0)[rb]{{\large
x}}} \put(85,27){\makebox(0,0)[cb]{{\large y}}}
\put(116,27){\makebox(0,0)[lb]{{\large z}}}
\put(85,9){\makebox(0,0)[ct]{horizontal}}
\put(15,9){\makebox(0,0)[ct]{vertical}}\end{picture}\end{center}

\begin{center}\begin{picture}(85,37)
\put(5,5){\dashbox{5}(10,30)[cc]{}} \put(20,5){\dashbox{5}(10,30)[cc]{}}
\put(50,3){\dashbox{5}(30,15)[cc]{}} \put(50,22){\dashbox{5}(30,15)[cc]{}}
\put(5,5){\line(1,1){10}} \put(5,15){\line(1,1){10}}
\put(20,15){\line(1,1){10}}
\put(20,25){\line(1,1){10}} \put(50,3){\line(2,1){30}}
\put(50,22){\line(2,1){30}}
\put(10,20){\circle*{1.5}} \put(10,10){\circle*{1.5}}
\put(25,20){\circle*{1.5}}
\put(25,30){\circle*{1.5}} \put(60,8){\circle*{1.5}}
\put(60,27){\circle*{1.5}}
\put(70,32){\circle*{1.5}} \put(70,13){\circle*{1.5}}
\put(40,20){\makebox(0,0)[cc]{$=$}} \put(43,2){\makebox(0,0)[ct]{Ehresmann
Axiom:
`middle four' exchange.}}
\end{picture}\end{center}

\begin{def8}[Modification, B\'enabou 1967] Let $s,t\in G_1(f,g).$ A
modification $\mu\in G_2(s,t)$ is a collection $(i+2)$-cells $\mu
x\in F_{i+1}(sx,tx),$ intertwining $s$ with $t,$ such that
$\forall\,a,b\in E_i(x,y)\subset E_{i+1}$ and
$$\forall\,\alpha\in E_i(x,y)(a,b)\subset E_{i+1}(a,b)\subset E_{i+2},$$
the following equations hold $\forall\,i\in\N,$

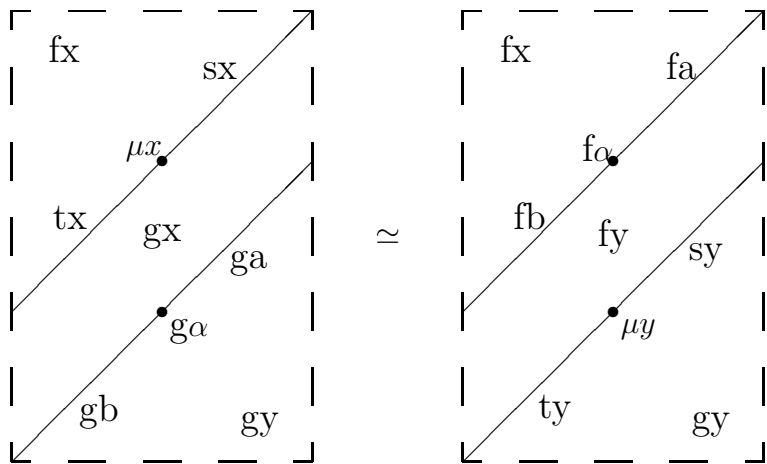
\begin{center}\begin{picture}(110,70)
\put(25,45){\circle*{1.5}} \put(25,25){\circle*{1.5}}
\put(85,45){\circle*{1.5}}
\put(85,25){\circle*{1.5}} \put(12,60){\makebox(0,0)[cc]{{\large fx}}}
\put(72,60){\makebox(0,0)[cc]{{\large fx}}}
\put(38,10){\makebox(0,0)[cc]{{\large
gy}}} \put(98,10){\makebox(0,0)[cc]{{\large gy}}}
\put(25,35){\makebox(0,0)[cc]{{\large gx}}}
\put(55,35){\makebox(0,0)[cc]{$\simeq$}}
\put(85,35){\makebox(0,0)[cc]{{\large fy}}}
\put(35,56){\makebox(0,0)[rb]{{\large
sx}}} \put(96,56){\makebox(0,0)[rb]{{\large fa}}}
\put(15,36){\makebox(0,0)[rb]{{\large tx}}}
\put(34,33){\makebox(0,0)[lt]{{\large
ga}}} \put(76,36){\makebox(0,0)[rb]{{\large fb}}}
\put(95,34){\makebox(0,0)[lt]{{\large sy}}}
\put(14,14){\makebox(0,0)[lt]{{\large
gb}}} \put(75,14){\makebox(0,0)[lt]{{\large ty}}}
\put(25,45){\makebox(0,0)[rb]{$\mu
x$}} \put(85,45){\makebox(0,0)[rb]{{\large f}$\alpha$}}
\put(26,24){\makebox(0,0)[lt]{{\large g}$\alpha$}}
\put(86,24){\makebox(0,0)[lt]{$\mu
y$}} \put(5,5){\dashbox{5}(40,60)[cc]{}}
\put(65,5){\dashbox{5}(40,60)[cc]{}}
\put(5,5){\line(1,1){40}} \put(65,5){\line(1,1){40}}
\put(5,25){\line(1,1){40}}
\put(65,25){\line(1,1){40}}\end{picture}\end{center}

\begin{gather*}x,y\in E_i,\quad a,b\in E_{i+1},\quad\alpha\in E_{i+2},\\
\begin{alignat*}{2}\text{surfaces:}&\quad\qquad fx,\;gx,\;fy,\;gy&&\in
F_i,\\
\text{arrows:}&\quad sx,tx,sy,ty,fa,fb,ga,gb&&\in F_{i+1},\\
\text{nodes:}&\quad\qquad\mu x,\;\mu y,\;f\alpha,\;g\alpha&&\in F_{i+2}.
\end{alignat*}\end{gather*}\end{def8}

Each of the above graphs, two elements of an operad, contains four paths
from $fx$ to $gy,$ bypassing over or under nodes. The above relation in
an operad tells that these four paths must be equal,
\begin{gather*}(ga)\circ(sx)=(sy)\circ(fa),\\
(gb)\circ(sx)=(sy)\circ(fb),\\(ga)\circ(tx)=(ty)\circ(fa),\\
(gb)\circ(tx)=(ty)\circ(fb),\\
\text{and}\quad``(g\alpha)\circ(\mu x)=(\mu y)\circ(f\alpha)''.
\end{gather*}
Therefore Definition 2.4 include Definition 2.3. Definition 2.4 is for
absolutely covariant $f$ and $g.$

Let $\mu,\nu\in G_0(c_E,c_F)(f,g)(s,t)\subset G_1(f,g)(s,t)\subset
G_2(s,t)\subset G_3.$

A `modification of modifications' $A\in G_3(\mu,\nu)\subset G_4$ is a
collection of $(i+3)$-cells $Ax\in F_{i+2}(\mu x,\nu x),$ intertwining
$\mu$ with $\nu$ such that the set of relations hold $\forall\,i\in\N.$
These relations needs graphics with 3-dimensional volumes, this time
$x,y\in E_i$ need to be represented by 3-dimensional volumes and $Ax$
and $Ay$ by nodes.

\section{n-Categories admissible by n-graph}
\begin{prob1} Given n-graph $G.$ How many (non isomorphic) n-category
structures are allowed?\end{prob1}

\begin{exm}
Let $G$ be a 2-graph whose 0-cells are sets, 1-cells are mappings between
sets, 2-cells are changes of mappings with the same domain and codomain.
Then every 3-cell is trivial.\end{exm}

\begin{exm}A 2-graph of 1-categories is a 2-category whose 1-cells are
functors and 2-cell are natural transformations.\end{exm}

\begin{exm} Let $G$ be a 1-graph such that every 0-cell $x\in G_{0}$ is
an (n-1)-category, and whose i-cells are (n-1)-category morphisms. If
there is an (associative) composition
\begin{equation}\begin{array}{c}
c=\{G_{0}(x,y)\times G_{0}(y,z)\overset{c_{x,y,z}}{\longrightarrow }%
G_{0}(x,z)\},
\end{array}
\end{equation}
then we obtain an n-category $G^{\prime }$ with only one 0-cell $G$ and
whose i-cells are (i-1)-cell of $G$.
\end{exm}

\begin{dfn}[Cobordism of manifolds, Stong 1998]
A manifold $M$ is said to be a cobordism from $A$ to $B$ if exists a
diffeomorphism from a disjoint sum, $\varphi \in \text{diff}(A^\ast\cup
B,\partial M).$ Two cobordisms $M(\varphi )$ and $M^{\prime }(\varphi
^{\prime })$ are equivalent if there is a
$\Phi\in\text{diff}(M,M^{\prime})$ such that $\varphi ^{\prime }=\Phi
\circ \varphi .$ The equivalence class of cobordisms is denoted by
$M(A,B)\in Cob(A,B).$\end{dfn}

\begin{exm} The unit cobordism $\partial (\Sigma \times \lbrack
0,1])=\Sigma \cup \Sigma^{\ast }.$\end{exm}

Nontrivial examples of cobordism are given by Ionicioiu [1997].

The 1-graph of oriented cobordisms $Cob$ has oriented manifolds of the
fixed dimension as objects, $Cob_0,$ and $Cob_1$ are cobordisms. A
$Cob_0$ is a Grothendieck rig $\equiv\;\mathbb{N}$-algebra,
\textit{i.e.}\, an commutative monoid with duality, a monoidal structure
is given by direct union `$\cup$' and the duality is reversal of
orientation.

Composition of cobordisms $c_{Cob}$ comes from gluing of manifolds. Let
$\varphi^{\prime}\in\text{diff}(C^*\cup D,\partial N).$ One can glue
cobordism $M$ with $N$ by identifying $B$ with $C^*,$ $(\varphi^{%
\prime})^{-1}\circ\varphi\in\text{diff}(B,C^*).$ We obtain the glued
cobordism $(M\circ N)(A,D)$ and a semigroup operation,
\begin{equation*}
c(A,B,D):Cob(A,B)\times Cob(B,D)\longrightarrow Cob(A,D).
\end{equation*}

\begin{exm} The 1-graph of oriented cobordisms $Cob$ with the
composition $c_{Cob}$ is a 1-category.\end{exm}

A \textit{surgery} is an operation of cutting a manifold $M$ and gluing
to cylinders. A surgery gives new cobordism: from $M(A,B)$ into
$N(A,B).$ The disjoint sum of $M(A,B)$ with $N(C,D)$ is a cobordism
$(M\cup N)(A\cup C,B\cup D).$ We got a 2-graph of cobordism $Cob$ with
$Cob_0=Man_d,$ $Cob_1=Man_{d+1},$ whose 2-cells from $Cob_2$ are surgery
operations.

Kerler [1998] found examples of categories formed by classes of cobordism
manifolds preserving disjoint sum or surgery. These examples are discussed
by Baez and Dolan [1996].

\begin{exm} Let $G$ be an $2$-graph such that every 0-cell is a k-graph,
every
1-cell $f\in G_{0}(x,y)$ is a k-graph morphism, composition of 1-cells is
defined as
a composition of k-graph morphisms. If every 2-cell $t\in G_{1}(f,g)$ is a
k-graph
transformation intertwining $f$ with $g$, where $f,g\in G_{0}(x,y)$, and
every k-cell$
k=3,\ldots ,$ is trivial, then the $2$-graph $G$ becomes a
2-category.\end{exm}

\begin{exm} Let $G$ be an $3$-graph such that every $x\in G_0$ be a
k-category
with a composition $c_x.$ Let $f,g\in G_{0}(c_x,c_y)$ be a k-category
morphisms, all
2-cells $s,t\in G_{1}(f,g)$ be a k-category transformations, and let every
3-cell
$\mu\in G_{2}\left(s,t\right)$ be a modification. Compositions of k-category
morphisms, transformations and modifications define compositions making the
3-graph
$G$ a 3-category.\end{exm}

\end{document}